# A scaling law for numerical computation of fracture stress intensity factors

Ronghao Cui[a,*], Jan Martin Nordbotten[a,b], Ivar Stefansson[a], Inga Berre[a]

[a] Center for Modeling of Coupled Subsurface Dynamics, Department of Mathematics, University of Bergen, Bergen, Norway

[b] NORCE Norwegian Research Center, Bergen, Norway

[*] Corresponding author. Email address: ronghao.cui@uib.no (R. Cui).

## Highlights

- Nearest-tip sampling does not ensure good convergence of stress intensity factors.
- A scaling law links sampling to mesh size for convergent stress intensity factors.
- Numerical benchmarks of mode I, mode II, and mixed modes validate the approach.
- Convergence analysis is performed for different numerical methods.

## Abstract

Numerical simulations of fracture propagation require accurate computation of stress intensity factors (SIFs) to determine whether a pre-existing fracture will propagate. Even though displacement solutions from well-established numerical methods generally converge with mesh refinement, the convergence behavior of SIF estimates remains less understood. In this work, we propose a scaling law to improve the performance of SIF extraction using two-point displacement correlation methods. This scaling law relates sampling distance from the fracture tip to mesh resolution. It is developed by balancing the asymptotic truncation error against the numerical discretization error. Three numerical discretization methods—standard quadratic finite elements, quadratic finite elements with quarter-point elements, and a multi-point stress approximation finite volume method—are employed to validate the proposed scaling law on fracture benchmark problems under pure-mode and mixed-mode loading. We find that the conventional nearest-tip sampling strategy fails to compute convergent SIF estimates. In contrast, the proposed scaling-based approach yields decreasing errors in the SIF estimates as the mesh is refined. Across all three discretization methods, the proposed scaling law captures the sampling locations required for convergent SIF extraction, suggesting a general numerical characteristic of near-tip SIF extraction.

## 1. Introduction

Fracture propagation arises in a broad range of engineering and natural applications, including subsurface energy systems, geotechnical engineering, structural materials, and natural hazards [1-5]. Accurate simulations of fracture propagation in deformable fractured media hinge on physics-based fracture criteria [6-8]. In linear elastic fracture mechanics, fracture propagation is understood from an energy perspective: the fracture starts to grow when the elastic strain energy released during fracture propagation equals or exceeds the surface energy required to extend the fracture [9]. The rate of change in the elastic strain energy with respect to the fracture surface area, defined as the energy release rate, is commonly used to characterize fracture propagation. The energy release rate is related to the stress intensity factor (SIF) [10]. In practical modeling of fracture mechanics, SIFs are key parameters for determining whether a pre-existing fracture propagates. Three different fracture modes are distinguished: opening (mode I), sliding (mode II), and tearing (mode III). The corresponding SIFs are defined as amplitudes of the leading $1/\sqrt{r}$-order terms in the asymptotic stress field near the fracture tip, where $r$ is the radial distance from the fracture tip [10, 11].Accurate characterization of near-tip fields is essential for evaluating SIFs, but remains challenging numerically. In particular, the stress field near the tip exhibits singular behavior as $r$ approaches zero. Connected with the stress field through the linear elastic constitutive relation and the strain-displacement relation, the displacement field has a consistent asymptotic structure, but exhibits $\sqrt{r}$-type behavior [12].

In numerical formulations of linear elastic fracture mechanics, SIFs are extracted through a post-processing procedure based on computed displacement, stress, or energy fields in the vicinity of the fracture tip. Two major approaches to SIF extraction, among others, are asymptotic-field methods and energy-based methods. On the one hand, asymptotic-field methods include displacement correlation and stress correlation methods, both of which depend on asymptotic expansions of the displacement and stress fields near the fracture tip in linear elastic fracture mechanics [13, 14]. On the other hand, energy-based methods leverage the path-independent $J$-integral, evaluated along a closed contour surrounding the fracture tip [15] or over an equivalent domain [16]. Extensions of the $J$-integral, commonly referred to as interaction integrals, are used in mixed-mode fracture problems to separate the contributions from modes I, II, and III and thereby compute the corresponding SIFs [17]. Compared with energy-based methods and stress correlation methods, displacement correlation methods are computationally efficient for simulations of fracture propagation and its coupling with fluid flow or heat transfer [18, 19].

Fracture propagation with SIF-based criteria has been extensively investigated in the family of finite element methods (FEMs) [20-23]. Classical FEM-based extensions to improve the computation of SIFs include enrichment functions and quarter-point elements (QPEs). Enrichment functions for displacement approximation are incorporated into extended finite element methods and allow the fracture to pass through non-conforming meshes without remeshing [24-26]. QPEs are introduced into the quadratic finite element method (QFEM) in which the mesh conforms to the fracture surfaces [22]. This modified formulation of the standard QFEM is referred to as QFEM-QPE. QFEM-QPE embeds the theoretical singularity at nearest-tip elements and improves the near-tip approximation by shifting midside nodes connected to the fracture tip [27, 28]. Recent developments have also made finite volume methods (FVMs) increasingly appealing in the context of fracture propagation [29-31]. As naturally conservative methods, FVMs are now used in coupled modeling of deformation, fluid flow, and heat transfer processes [29-32]. In particular, the multi-point stress approximation (MPSA) was introduced to show displacement accuracy comparable to that of the lowest-order finite element methods [33]. To improve the discretization robustness on simplex grids, the weakly symmetric variant of MPSA was subsequently proposed [34]. Nevertheless, accurate extraction of SIF estimates from FVM formulations such as MPSA has not yet been well established.

Beyond method-specific accuracy of SIF extraction, a more fundamental issue across different numerical methods is the convergence of SIF estimates with mesh refinement. The simplest displacement correlation approach to extracting SIFs uses a single displacement jump solution [35]. We refer to this type of method as the single-point method. Previous studies found that a linear correction to the single-point method could effectively reduce the error in SIF estimates by using two samples of the displacement jump along the fracture [27, 35]. One standard strategy is to sample nearest-tip nodes, as done in QFEM-QPE [14, 20, 22]. Displacement solutions obtained from well-established numerical methods, such as FEMs and FVMs, generally converge with mesh refinement. However, the convergence behavior of SIF estimates extracted from numerical solutions near the fracture tip is often less straightforward [13, 14, 36]. The major source of this discrepancy across different numerical methods lies in the SIF extraction procedure itself, owing to the low regularity of the solution near the fracture tip [13, 14, 36]. Even though the multiple-point displacement correlation methods proposed by Chan et al. [13] can partially alleviate the sensitivity of SIF estimates to numerical errors by avoiding the direct use of the nearest-tip samples, it remains unclear how sampling locations should be determined for convergent SIF estimates.

In this work, we propose a scaling law that relates sampling distance to mesh resolution for addressing the common challenge of SIF estimation across three representative numerical discretization methods: QFEM, QFEM-QPE, and MPSA. Section 2 introduces the general problem statement of linear elastic fracture mechanics and the formulation of SIF extraction. Section 3 presents the theoretical motivation and mathematical formulation of the proposed scaling law. Section 4 describes the discretization details of QFEM, QFEM-QPE, and MPSA for two-dimensional (2D) benchmark problems involving mode I and mode II. These benchmark problems allow the proposed scaling law to be evaluated in pure modes and mixed modes. Section 5 investigates the limitations of the nearest-tip sampling strategy and evaluates the performance of the proposed scaling law across the selected discretization methods. The results show that the sampling locations that yield convergent SIF estimates can be captured by the proposed scaling law, revealing a general numerical characteristic of near-tip SIF extraction.

## 2. Linear elastic fracture mechanics formulation

### 2.1 Governing equations

We consider a two-dimensional, homogeneous, and linear elastic body occupying the domain $\Omega$, which represents a deformable fractured medium (see Fig. 1). With the constitutive relation governed by Hooke's law, the momentum balance equation for linearized elasticity reads

$$\begin{aligned} &\nabla \cdot \boldsymbol{\sigma} + \boldsymbol{f} = 0 && \text{in } \Omega, \\ &\boldsymbol{\sigma} = \mathbb{C} : \tfrac{1}{2}\left(\nabla \boldsymbol{u} + (\nabla \boldsymbol{u})^{\mathrm{T}}\right) && \text{in } \Omega, \end{aligned} \tag{1}$$

where $\boldsymbol{f}$ is the body force, $\boldsymbol{\sigma}$ is the stress tensor, $\mathbb{C}$ is the stiffness tensor, and $\boldsymbol{u}$ is the displacement field. For isotropic materials, the constitutive relation in Eq. (1) reduces to $\boldsymbol{\sigma} = \lambda(\nabla \cdot \boldsymbol{u})\boldsymbol{I} + \mu\left(\nabla \boldsymbol{u} + (\nabla \boldsymbol{u})^{\mathrm{T}}\right)$, where $\lambda$ and $\mu$ are Lamé parameters. The external boundary of the domain, denoted by $\Gamma_D$, is subject to the Dirichlet boundary condition given by

$$\boldsymbol{u} = \boldsymbol{u}_D \quad \text{on } \Gamma_D, \tag{2}$$

where $\boldsymbol{u}_D$ is the prescribed boundary displacement. The two opposing surfaces of a fracture are treated as internal traction boundaries of the domain, $\Gamma_s^+$ and $\Gamma_s^-$. The tractions imposed on $\Gamma_s^+$ and $\Gamma_s^-$, $\boldsymbol{t}_s^+$ and $\boldsymbol{t}_s^-$, are given by

$$\begin{aligned} \boldsymbol{\sigma} \cdot \boldsymbol{n}_s^+ &= \boldsymbol{t}_s^+ \quad \text{on } \Gamma_s^+, \\ \boldsymbol{\sigma} \cdot \boldsymbol{n}_s^- &= \boldsymbol{t}_s^- \quad \text{on } \Gamma_s^-, \end{aligned} \tag{3}$$

where $\boldsymbol{n}_s^+$ and $\boldsymbol{n}_s^-$ are outward unit normal vectors on $\Gamma_s^+$ and $\Gamma_s^-$, respectively. For traction-free fractures, $\boldsymbol{t}_s^+ = \boldsymbol{t}_s^- = 0$. If the fracture is subject to a constant pressure $p_0$, the pressure is considered as the constant traction acting normal to the fracture surfaces, i.e., $\boldsymbol{t}_s^+ = -p_0 \boldsymbol{n}_s^+$ and $\boldsymbol{t}_s^- = -p_0 \boldsymbol{n}_s^-$. The displacement jump across the two sides of the fracture is defined as

$$[\![\boldsymbol{u}]\!] = \boldsymbol{u}_s^+ - \boldsymbol{u}_s^-, \tag{4}$$

where $\boldsymbol{u}_s^+$ and $\boldsymbol{u}_s^-$ are displacements on $\Gamma_s^+$ and $\Gamma_s^-$, respectively. The jump across the fracture, $[\![\boldsymbol{u}]\!]$, consists of normal and tangential components, $[\![\boldsymbol{u}]\!]_n$ and $[\![\boldsymbol{u}]\!]_t$. The above formulation defines the classical near-tip problem in linear elastic fracture mechanics for traction-free fractures and uniformly pressurized fractures. The geometry for this problem is illustrated in Fig. 1.

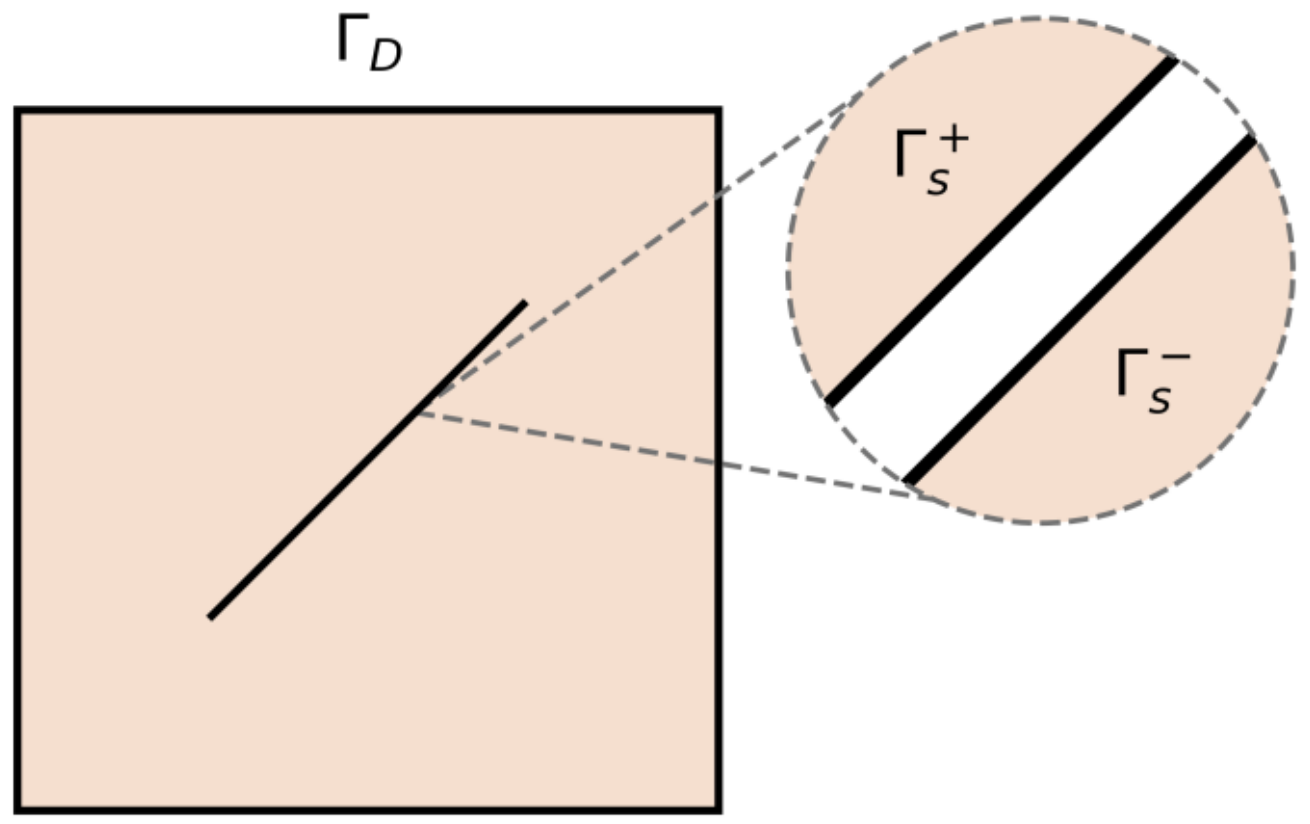

**Fig. 1.** The fractured domain with given external and internal boundaries for the near-tip problem.

### 2.2 Near-tip SIF extraction

The SIFs are defined through the leading-order terms in the asymptotic expansions of the stress or displacement fields near a fracture tip [11]. In a local polar coordinate system centered at the fracture tip, the two opposing surfaces of a straight fracture correspond to $\theta = \pm\pi$, where $\theta$ is the angle measured from the fracture plane. This simplifies asymptotic expansions of the displacement jump across the two fracture surfaces [14], yielding

$$\begin{aligned} [\![\boldsymbol{u}(r)]\!]_n &= \frac{(\kappa+1)}{\mu} K_I \sqrt{\frac{r}{2\pi}} + \sum_{m=1}^{\infty} f_m^I \, r^{m+1/2}, \\ [\![\boldsymbol{u}(r)]\!]_t &= \frac{(\kappa+1)}{\mu} K_{II} \sqrt{\frac{r}{2\pi}} + \sum_{m=1}^{\infty} f_m^{II} \, r^{m+1/2}, \end{aligned} \tag{5}$$

where $K_I$ and $K_{II}$ are exact SIFs of mode I and mode II, respectively, $f_m^I$ and $f_m^{II}$ are coefficients of higher-order terms in mode I and mode II, respectively, and $\kappa$ is the Kolosov constant. The value of $\kappa$ depends on the assumption of plane strain ($\kappa = 3 - 4\nu$) or plane stress ($\kappa = (3-\nu)/(1+\nu)$), where $\nu$ is Poisson's ratio. It should be noted that the elastic properties in our two-dimensional formulation may be parameterized by Lamé parameters $\mu$ and $\lambda$, or equivalently by Young's modulus $E$ and Poisson's ratio $\nu$. However, the conversion between these two sets of parameters depends on the assumption of plane strain or plane stress. The exact SIFs can be given asymptotically in terms of displacement jumps evaluated at a given radial distance from the fracture tip $r$ with corresponding higher-order corrections

$$K_I = \frac{\mu}{(\kappa+1)}\sqrt{\frac{2\pi}{r}}[\![\boldsymbol{u}(r)]\!]_n + \sum_{m=1}^{M} g_m^I r^m + O(r^{M+1}),$$
$$K_{II} = \frac{\mu}{(\kappa+1)}\sqrt{\frac{2\pi}{r}}[\![\boldsymbol{u}(r)]\!]_t + \sum_{m=1}^{M} g_m^{II} r^m + O(r^{M+1}), \tag{6}$$

where $M \geq 1$ is the truncation order of the asymptotic expansion and

$$g_m^I = -\frac{\mu}{(\kappa+1)}\sqrt{2\pi} f_m^I,$$
$$g_m^{II} = -\frac{\mu}{(\kappa+1)}\sqrt{2\pi} f_m^{II}. \tag{7}$$

The first higher-order term in Eq. (6) after the leading term results in an $O(r)$ truncation error in the extracted mode-I and mode-II SIFs. These relations hold for both pure-mode and mixed-mode problems, providing the asymptotic basis for displacement correlation methods. Let $r_i$ denote the radial distance of sampling point $i$ from the fracture tip. SIFs can be approximated by the single-point displacement correlation method (denoted by the superscript *sp*), as given by

$$K_I^{sp}(r_i) = \frac{\mu}{(\kappa+1)}\sqrt{\frac{2\pi}{r_i}}[\![\boldsymbol{u}(r_i)]\!]_n,$$
$$K_{II}^{sp}(r_i) = \frac{\mu}{(\kappa+1)}\sqrt{\frac{2\pi}{r_i}}[\![\boldsymbol{u}(r_i)]\!]_t. \tag{8}$$

Assuming that the leading term in Eq. (6) correctly accounts for the singularity, the accuracy of the SIF estimates can be improved by keeping terms up to the first order. Using two distinct sampling points $r_i$ and $r_j$ along the same radial direction, the linear term can be eliminated. The resulting approximation of the SIFs using a two-point displacement correlation method, denoted by the superscript *tp*, is given by

$$K_I^{tp}(r_i, r_j) = K_I^{sp}(r_i) + \frac{r_i}{r_j - r_i}\left(K_I^{sp}(r_i) - K_I^{sp}(r_j)\right),$$
$$K_{II}^{tp}(r_i, r_j) = K_{II}^{sp}(r_i) + \frac{r_i}{r_j - r_i}\left(K_{II}^{sp}(r_i) - K_{II}^{sp}(r_j)\right). \tag{9}$$

## 3. Formulation of the scaling law

In this section, we develop a heuristic argument for selecting the sampling distance for SIF extraction using displacement correlation methods as the mesh is refined. SIF estimation is subject to two sources of error: asymptotic truncation error and numerical discretization error. As a result, the accuracy of displacement-jump-based SIF estimation is governed by two characteristic length scales: the mesh size $h$ and the sampling distance from the tip $r$. Specifically, when $r$ is large, the asymptotic truncation error dominates. When $r$ is small, the local numerical solution is strongly affected by the discretization error due to the fracture-tip singularity. Therefore, for a given mesh resolution, it is reasonable to assume that an optimal sampling region exists where the total error arising from these two sources is minimized.

The fracture-tip asymptotic expansions in Eq. (5) are analytically derived from the governing equation and the geometric discontinuity in problems of linear elastic fracture mechanics. Eq. (5) indicates that the leading-order term of the displacement-jump expansion is proportional to $\sqrt{r}$. However, the neglected higher-order terms introduce the asymptotic truncation error. Using the SIF estimates by the two-point displacement correlation method, the asymptotic truncation error $e_{\varphi}^{asym}$ for the fracture mode $\varphi \in \{I, II\}$ scales as

$$e_{\varphi}^{asym} = O(r^2). \tag{10}$$

A detailed explanation for Eq. (10) is included in Appendix A.

In addition, the numerical estimation of the displacement jump is affected by method-dependent discretization error. We understand this as follows. For a finite element method based on $p$th-order polynomials, the global discretization error for a sufficiently smooth solution in the domain can be bounded by Céa's lemma as

$$||\nabla(\boldsymbol{u} - \boldsymbol{u}_h)||_{L^2(\Omega)} \leq ch^q ||\nabla^{q+1}\boldsymbol{u}||_{L^2(\Omega)}, \quad q \leq p, \tag{11}$$

where $q$ is the convergence order permitted by both the polynomial degree and the Sobolev regularity of the exact solution, $\boldsymbol{u}_h$ is the finite element displacement field, and $c$ is a positive constant independent of $h$. At the fracture tip (i.e., $r = 0$), the derivatives of $\boldsymbol{u}$ are not bounded for $q \geq \frac{1}{2}$; thus the highest convergence rate that can be expected from any finite element method is (slightly less than) $O(h^{1/2})$ for the gradient approximation. This limitation is a best approximation property, and can also be considered a general expectation for other classes of numerical methods.

Nevertheless, away from the singularity, optimal-order convergence rates (i.e. $q = p$) have been observed in practical computations [37]. Since the stress intensity factors in Eq. (6) are evaluated from the displacement jump instead of the displacement gradient, we now consider the approximation error of the displacement field. For a sufficiently smooth solution, a finite element method based on $p$th-order polynomials typically yields an $L^2$-norm error in the displacement field of order $O(h^{p+1})$. Since the displacement jump is obtained from the traces of the displacement field on the two fracture surfaces, we are interested in evaluating the displacement solution on the boundary. We expect the discretization error of the displacement jump to scale as $O(h^{p+1-\delta})$, where $\delta$ accounts for the loss of convergence order resulting from evaluating displacement traces on the fracture surfaces. When sufficiently far away from the fracture tip, its nominal values should be expected to be around $1/2$.

Since the global $L^2$-norm error estimate does not by itself determine the errors at individual sampling points, we need an additional assumption for local errors. The signed discretization errors in normal and tangential components of displacement jump, $\eta_I^{num}(r)$ and $\eta_{II}^{num}(r)$, are respectively defined as

$$\begin{aligned} \eta_I^{num}(r) &= [\![\boldsymbol{u}_h(r)]\!]_n - [\![\boldsymbol{u}(r)]\!]_n, \\ \eta_{II}^{num}(r) &= [\![\boldsymbol{u}_h(r)]\!]_t - [\![\boldsymbol{u}(r)]\!]_t. \end{aligned} \tag{12}$$

To obtain a scaling estimate for the discretization error, we assume that the displacement-jump errors at sampling locations, $r_i$ and $r_j$, have comparable magnitudes, sharing the same leading-order coefficient. Then, $\eta_\varphi^{num}(r_i)$ and $\eta_\varphi^{num}(r_j)$ for the fracture mode $\varphi \in \{I, II\}$ are given by

$$\eta_\varphi^{num}(r_i) = c_\varphi\, h^s, \eta_\varphi^{num}(r_j) = c_\varphi\, h^s, \tag{13}$$

where $s$ is the local convergence exponent of the sampled displacement jump and $c_\varphi$ is the signed coefficient. The convergence exponent $s$ in Eq. (13) assumes the sampled errors inherit the convergence behavior of the displacement jump field without local error variance. Under this assumption, the analysis in Appendix B indicates that the numerical discretization error $e_\varphi^{num}$ of the two-point displacement correlation method scales as

$$e_\varphi^{num} = O(h^s r^{-0.5}). \tag{14}$$

The total error in the SIF estimation, $e_\varphi^{tot} = e_\varphi^{asym} + e_\varphi^{num}$, is therefore assumed to scale as

$$e_\varphi^{tot} = O(r^2) + O(h^s r^{-0.5}). \tag{15}$$

Following the principle of error balancing, the optimal sampling distance is obtained by requiring the discretization error and the asymptotic truncation error to be of the same order

$$O(h^s r^{-0.5}) = O(r^2). \tag{16}$$

Based on this error-balancing condition in Eq. (16), we propose a scaling law in a dimensionless form for sampling the displacement jump near the fracture tip, as given by

$$\frac{r}{a} = C\left(\frac{h}{a}\right)^{\alpha}, \tag{17}$$

where $a$ is the half length of the fracture, $C$ is the scaling constant, and $\alpha$ is the scaling exponent, $\alpha = 2s/5$.

## 4. Numerical discretization and software

For the fracture problem stated in Section 2, we select three different discretization methods, QFEM, QFEM-QPE, and MPSA, for numerical simulations and convergence analysis. QFEM-QPE is the primary reference method because of its widespread use in SIF evaluation in fracture propagation simulations. In order to understand the impact of QPEs on the representation of the fracture-tip singularity, the standard QFEM is considered as a baseline finite element discretization method for SIF extraction. In addition, MPSA is included to investigate whether the observed numerical behavior is still valid for finite volume methods.

### 4.1 Finite element discretization

#### 4.1.1 Standard quadratic finite elements

The standard QFEM discretizes the linear elastic problem defined in Section 2 using a fracture-conforming simplex partition of the 2D matrix domain. The fracture surfaces are explicitly represented by internal boundaries on element edges, as illustrated in Fig. 2a. QFEM approximates the displacement field by quadratic nodal interpolation over six-node triangular elements. We denote the displacement approximation in quadratic finite elements by $\boldsymbol{u}_h^e$. The weak form of Eq. (1) in the FEM formulation is obtained by multiplying a test function $\boldsymbol{v}_h^e$ and integrating over the whole domain $\Omega$ by parts. Let $V_s$ be the solution space satisfying the prescribed Dirichlet boundary conditions on $\Gamma_D$ and let $V_t$ denote the corresponding test space. The weak form reads: find $\boldsymbol{u}_h^e \in V_s$ such that, for all $\boldsymbol{v}_h^e \in V_t$,

$$\int_{\Omega} \boldsymbol{\varepsilon}(\boldsymbol{v}_h^e) : \mathbb{C} : \boldsymbol{\varepsilon}(\boldsymbol{u}_h^e)\, dV = \int_{\Omega} \boldsymbol{v}_h^e \cdot \boldsymbol{f}\, dV + \int_{\Gamma_s^+} \boldsymbol{v}_h^e \cdot \boldsymbol{t}_s^+\, d\Gamma + \int_{\Gamma_s^-} \boldsymbol{v}_h^e \cdot \boldsymbol{t}_s^-\, d\Gamma, \tag{18}$$

where $\boldsymbol{\varepsilon}(\cdot)$ is the infinitesimal strain operator defined as

$$\boldsymbol{\varepsilon}(\boldsymbol{w}) = \frac{1}{2}\left(\nabla \boldsymbol{w} + (\nabla \boldsymbol{w})^{\mathrm{T}}\right), \tag{19}$$

for any vector-valued function $\boldsymbol{w}$. Over each six-node triangular element, $\boldsymbol{u}_h^e$ is interpolated in the reference coordinates $(\xi, \eta)$ as

$$\boldsymbol{u}_h^e = \sum_{i=1}^{6} N_i^{tri}(\xi, \eta)\boldsymbol{u}_i^e, \tag{20}$$

where $N_i^{tri}(\xi, \eta)$ is the quadratic shape function of node $i$ in six-node triangular elements, and $\boldsymbol{u}_i^e$ is the displacement vector at node $i$. The physical geometry of each element is defined by isoparametric mapping, as given by

$$\boldsymbol{x}_h^e(\xi, \eta) = \sum_{i=1}^{6} N_i^{tri}(\xi, \eta)\boldsymbol{x}_i^e, \tag{21}$$

where $\boldsymbol{x}_i^e$ is the nodal coordinate vector and $\boldsymbol{x}_h^e(\xi, \eta)$ is the corresponding position vector within the physical element. Substituting Eqs. (20) and (21) into the weak form in Eq. (18), we arrive at the assembled linear system

$$\boldsymbol{K}\boldsymbol{U}^e = \boldsymbol{F}, \tag{22}$$

where $\boldsymbol{K}$ is the global stiffness matrix, $\boldsymbol{U}^e$ is the vector of global nodal displacement unknowns, and $\boldsymbol{F}$ is the global load vector including contributions from body force and prescribed traction. Since the fracture surfaces are represented by two coincident internal boundaries with duplicated nodes, QFEM calculates the displacement jump across the fracture through corresponding nodal pairs.

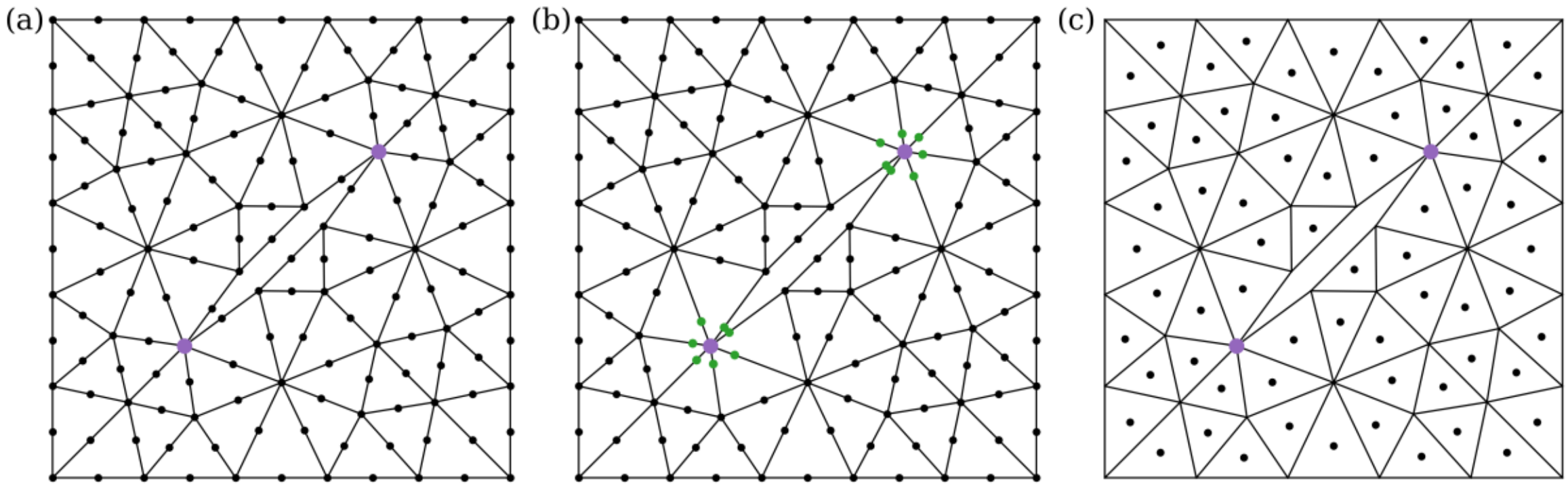

**Fig. 2.** Discretization illustrations of the 2D fractured domain by three numerical methods: (a) QFEM, (b) QFEM-QPE, and (c) MPSA. Black dots, purple dots, and green dots denote locations of displacement unknowns, fracture tips, and quarter-point nodes, respectively. Fracture tips coincide with displacement unknowns in both QFEM and QFEM-QPE, whereas they are not represented by displacement unknowns in MPSA. Quarter-point nodes in QFEM-QPE carry displacement degrees of freedom.

### 4.1.2 Quadratic finite elements with quarter-point elements

QFEM-QPE employs the same weak formulation and global finite element system as the standard QFEM. The only modification is made to elements adjacent to the fracture tips. For those fracture-tip elements, the midside nodes on the edges emanating from the fracture tip are shifted to the quarter-point locations measured from the tip, as shown in Fig. 2b. Each fracture-tip triangular element is constructed through a mapping from a collapsed eight-node quadrilateral element. The local nodes of the eight-node quadrilateral element are numbered consecutively along its boundary, with the odd-numbered nodes located at the corners and the even-numbered nodes at the midsides, starting from node $i = 1$ at one end of the edge to be collapsed. The physical coordinates of nodes $i = 1,7,8$ are mapped to the same fracture-tip location represented by node $i = 1$ [27]. The shape functions of eight-node quadrilateral elements are denoted as $N_i^{quad}(\xi, \eta)$, where $i = 1,2,3, \ldots ,8$. The isoparametric mappings of geometry and displacement are given by

$$\begin{aligned} \boldsymbol{u}_h^e &= \textstyle\sum_{i=1}^{6} N_i^*(\xi, \eta) \boldsymbol{u}_i^e, \\ \boldsymbol{x}_h^e &= \textstyle\sum_{i=1}^{6} N_i^*(\xi, \eta) \boldsymbol{x}_i^e, \end{aligned} \tag{23}$$

where $N_i^*(\xi, \eta)$ is the quadratic shape function of node $i$ within QPEs and is defined as

$$N_i^*(\xi, \eta) = \begin{cases} N_1^{quad}(\xi, \eta) + N_7^{quad}(\xi, \eta) + N_8^{quad}(\xi, \eta) \ \text{ when } \ i = 1, \\ N_i^{quad}(\xi, \eta) \ \text{ when } \ i = 2,3, \ldots ,6. \end{cases} \tag{24}$$

The geometric modification alters the isoparametric mapping within QPEs. Along an edge connected to the fracture tip, we denote the physical radial distance from the tip by $r$ and the corresponding reference distance as $\zeta$. In standard quadratic elements with midside nodes, the mapping between them is $r \propto \zeta$. However, the quarter-point modification yields the asymptotic mapping $\zeta \propto \sqrt{r}$. Consequently, the quadratically interpolated displacement $\boldsymbol{u}_h^e$ in terms of $r$ becomes

$$\boldsymbol{u}_h^e(r) = \boldsymbol{a}_0 + \boldsymbol{a}_1 \sqrt{r} + \boldsymbol{a}_2 r, \tag{25}$$

where $\boldsymbol{a}_0$, $\boldsymbol{a}_1$, and $\boldsymbol{a}_2$ are coefficient vectors. The collapsed eight-node quarter-point element embeds the square-root asymptotic displacement term into the approximation space. This enables QFEM-QPE to reproduce the classical near-tip behavior of linear elastic fracture mechanics with improved accuracy of near-tip solutions compared with the standard QFEM.

### 4.2 Finite volume discretization

Our finite volume discretization for this problem is carried out by the weakly symmetric variant of MPSA. The computational domain is represented by the same mesh as in Section 4.1, where each cell $i$ occupies a control volume $\Omega_i$. As shown in Fig. 2c, the cell-centered displacement, $\boldsymbol{u}_i^v$, is designated as the primary unknown [38]. The finite volume formulation integrates Eq. (1) over $\Omega_i$ and applies the divergence theorem [33, 34], yielding

$$\int_{\partial\Omega_i} \boldsymbol{\sigma} \cdot \boldsymbol{n}\, dA + \int_{\Omega_i} \boldsymbol{f}\, dV = 0. \tag{26}$$

where $\boldsymbol{n}$ is the outward normal vector on boundary $\partial\Omega_i$. The discrete momentum balance for cell $i$ is then written as

$$\sum_j \boldsymbol{T}_{i,j} + |\Omega_i| \boldsymbol{f}_i = 0, \tag{27}$$

where $\boldsymbol{T}_{i,j}$ is the integrated traction force acting on cell $i$ over the face shared by cell $i$ and cell $j$, $\boldsymbol{f}_i$ is the volume-averaged body force over cell $i$, and $|\Omega_i|$ is the volume of cell $i$. To construct the MPSA discretization, each cell is divided into vertex-based subcells. The vertex-based subcell of cell $i$ towards vertex $l$ is denoted by $\Omega_{i,l}$. The displacement field within $\Omega_{i,l}$ is approximated by a piecewise linear reconstruction, $\boldsymbol{u}_h^v$, based on the cell-centered displacement $\boldsymbol{u}_i^v$ and the local gradient $\boldsymbol{G}_{i,l}$, as given by

$$\boldsymbol{u}_h^v = \boldsymbol{u}_i^v + \boldsymbol{G}_{i,l}(\boldsymbol{x} - \boldsymbol{x}_i), \tag{28}$$

where $\boldsymbol{x}$ is the point within $\Omega_{i,l}$ and $\boldsymbol{x}_i$ is the center of cell $i$. The local gradients are determined by enforcing displacement continuity and traction continuity on subfaces shared by neighboring subcells. For the weakly symmetric MPSA formulation, the discrete stress in $\Omega_{i,l}$, the discrete stress tensor within $\Omega_{i,l}$, $\boldsymbol{\sigma}_{i,l}$, reads

$$\boldsymbol{\sigma}_{i,l} = \mathbb{C}_i : \boldsymbol{G}_{i,l} - \frac{\langle \mathbb{C}_i : \boldsymbol{G}_{i,l} \rangle_l - \langle \mathbb{C}_i : \boldsymbol{G}_{i,l} \rangle_l^{\mathrm{T}}}{2}, \tag{29}$$

where $\mathbb{C}_i$ is the stiffness tensor in $\Omega_i$ and $\langle \mathbb{C}_i : \boldsymbol{G}_{i,l} \rangle_l$ is the weighted average of discrete stress tensors over all subcells associated with vertex $l$. Solving the local linear system towards vertex $l$ eliminates $\boldsymbol{G}_{i,l}$ and yields a discrete relation in which $\boldsymbol{T}_{i,j}$ is a linear function of cell-centered displacement unknowns in the neighborhood. The coefficients of this linear relation are referred to as stress-weight tensors and assembled into the global linear system formulated by applying Eq. (27) to all cells.

The fracture surfaces are represented by internal boundary faces of matrix grids, where the prescribed traction on the fracture surfaces enters the global discretization through Neumann boundary conditions [39]. These local linear systems involving internal boundaries are constructed independently of the prescribed fracture loading. After the cell-centered displacement field is solved, the face-centered displacement values on the two surfaces of the fracture are reconstructed from the local linear approximation in Eq. (28), from which the displacement jump across the fracture is then calculated.

Existing convergence proofs for MPSA methods emphasize more robustness with respect to heterogeneous material parameters than extracting optimal convergence rates [40]. In contrast to this theoretical focus, numerical experiments on simplex grids have reported convergence rates of up to 2 for displacement and approximately 3/2 for stress [34].

### 4.3 The discretized scaling law

The proposed scaling formulation aims to provide a unified sampling strategy of displacement jump near the fracture tip for minimizing SIF estimation errors and making numerical SIF estimates converge across different discretization methods. Since QFEM, QFEM-QPE, and MPSA have different capabilities in characterizing the asymptotic fields near the fracture tip, the scaling parameters $C$ and $\alpha$ in Eq. (17) are not prescribed a priori but determined from numerical experiments. The displacement jump at a given sampling distance is calculated for all three discretization methods as the difference between paired displacement values at corresponding locations on the opposing fracture surfaces. On a known mesh, however, the sampling distance predicted by the scaling law generally does not coincide with the location of an available node or face center. Therefore, the closest available node or face center to the target distance calculated from the scaling law is used for SIF estimation in a given discretization method. For MPSA, the displacement unknowns on the fracture surfaces are reconstructed at paired face centers, which are chosen as sampling locations. By contrast, the sampling locations in QFEM and QFEM-QPE are node-wise because the displacement solutions are represented by nodal degrees of freedom. We employ the scaling law for the two-point displacement correlation methods, where only one sampling location is determined. The other sampling location is then chosen as the nearest neighbor closer to the fracture tip in a discrete manner.

### 4.4 Software

The simulation models in this work are carried out based on the open-source framework *PorePy* [41], while the QFEM and QFEM-QPE implementations are adapted from Dang et al. [42], with modifications for this work. All three discretization methods share the same simplex meshes generated within *PorePy* using Gmsh [43]. The source code for this work, together with the runscripts used to reproduce all numerical examples presented, is available in Cui et al. [44]. A generative artificial intelligence (AI) tool was used to provide limited support during code development. Specifically, the OpenAI-GPT5.6 model was used to assist with code formatting, readability, and debugging. All AI-suggested code modifications were carefully reviewed and verified to ensure compliance with scientific standards.

### 4.5 Benchmarks for SIF calculations

We consider four analytical benchmark problems for SIF calculations to cover mode-I, mode-II, and mixed-mode fracture behavior. They are Westergaard's mode-I problem, Westergaard's mode-II problem, Westergaard's mixed-mode problem, and Sneddon's mode-I problem, all of which provide analytical solutions of displacement jumps and SIFs to evaluate the convergence of the extracted SIF estimates. These comprise Westergaard-type traction-free fracture problems and the Sneddon pressurized fracture problem, as detailed in the following sections.

#### 4.5.1 Westergaard-type traction-free fracture problems

We adopt three Westergaard-type benchmarks, corresponding to mode I, mode II, and mixed-mode loading, respectively. In all three problems, the fracture surfaces are traction-free. Westergaard's mode-I and mode-II problems are both pure-mode fracture problems, both defined by a straight traction-free fracture in an infinite elastic medium. The mode-I and mode-II benchmarks have uniform far-field normal and in-plane shear stresses, $\sigma_\infty$ and $\tau_\infty$, respectively [45]. These two different loading conditions induce normal opening and tangential sliding between the two fracture surfaces, corresponding to SIFs of mode I and mode II at the fracture tip, respectively. For a fracture with half-length $a$, analytical solutions of normal and tangential displacement jumps, $[\![\boldsymbol{u}_{ref}^{w}]\!]_n$ and $[\![\boldsymbol{u}_{ref}^{w}]\!]_t$, along the fracture are respectively written as

$$[\![\boldsymbol{u}_{ref}^{w}]\!]_n = \frac{(\kappa+1)\sigma_\infty a}{2\mu}\sqrt{1-\frac{x_f^2}{a^2}}, \ [\![\boldsymbol{u}_{ref}^{w}]\!]_t = \frac{(\kappa+1)\tau_\infty a}{2\mu}\sqrt{1-\frac{x_f^2}{a^2}}, \tag{30}$$

where $x_f$ is the distance from the fracture center. The analytical SIF solutions for mode I and mode II, $K_{\mathrm{I},ref}^{w}$ and $K_{\mathrm{II},ref}^{w}$, are given by

$$K_{I,ref}^{w} = \sigma_\infty\sqrt{\pi a}, \ \ K_{II,ref}^{w} = \tau_\infty\sqrt{\pi a}. \tag{31}$$

If the fractured domain is subject to both $\sigma_\infty$ and $\tau_\infty$, the fracture undergoes simultaneously normal opening and tangential sliding, leading to a mixed-mode behavior with both $K_I$ and $K_{II}$ contributions. This benchmark is referred to as Westergaard's mixed-mode problem.

### 4.5.2 Sneddon pressurized fracture problem

For Sneddon's mode-I problem, a single fracture is subject to a constant pressure $p_0$ acting as a normal traction on the fracture surfaces [46]. This problem only leads to normal opening of the fracture, for which the analytical normal displacement jump $[\![\boldsymbol{u}_{ref}^{s}]\!]_n$ is

$$[\![\boldsymbol{u}_{ref}^{s}]\!]_n = \frac{(\kappa+1)p_0 a}{2\mu}\sqrt{1-\frac{x_f^2}{a^2}}. \tag{32}$$

The corresponding analytical solution of the SIF in mode I is

$$K_{I,ref}^{s} = p_0\sqrt{\pi a}. \tag{33}$$

### 4.5.3 Benchmark implementation based on boundary element methods

Westergaard's and Sneddon's analytical solutions are both derived in infinite domains, while numerical simulations in this work have to be implemented in finite computational domains. We employ the boundary element method (BEM) to construct the reference displacement fields for these benchmarks. This allows us to calculate reference displacements at arbitrary locations in the matrix domain. The reference displacement field $\boldsymbol{u}_{ref}$ is decomposed into two contributions, the homogeneous far-field displacement $\boldsymbol{u}_{\infty}$ and the fracture-induced perturbation displacement $\boldsymbol{u}_f$, i.e.

$$\boldsymbol{u}_{ref} = \boldsymbol{u}_{\infty} + \boldsymbol{u}_f. \tag{34}$$

During the benchmark implementation, the Dirichlet boundary displacement $\boldsymbol{u}_D$ in Eq. (2) is obtained by evaluating $\boldsymbol{u}_{ref}$ on $\Gamma_D$.

In Westergaard's problems, the far-field normal strain and shear strain, $\varepsilon_n$ and $\varepsilon_{nt}$, correspond to $\sigma_\infty$ and $\tau_\infty$, respectively, through linear elasticity

$$\varepsilon_n = \frac{\sigma_\infty}{2(\mu+\lambda)}, \ \varepsilon_{nt} = \frac{\tau_\infty}{2\mu}. \tag{35}$$

The background field for the far-field loading is calculated by the homogeneous linear elastic solution

$$\boldsymbol{u}_{\infty} = \varepsilon_n(r_t\boldsymbol{t} + r_n\boldsymbol{n}) + \varepsilon_{nt}(r_n\boldsymbol{t} + r_t\boldsymbol{n}), \tag{36}$$

where $r_n$ and $r_t$ are the normal and tangential coordinates measured in the local fracture coordinate system, respectively, and $\boldsymbol{n}$ and $\boldsymbol{t}$ are unit normal and unit tangential vectors of the fracture, respectively. For Sneddon's problem, $\boldsymbol{u}_{\infty} = 0$.

The perturbation contribution $\boldsymbol{u}_f$ induced by the fracture itself is computed by the displacement-discontinuity boundary element method [47]. The analytical displacement-jump distributions given in Eqs. (30) and (32) are discretized into a set of boundary-element segments along the fracture. Each segment is represented as a constant displacement-discontinuity source characterized by its normal opening and tangential sliding. The perturbation contribution generated by each segment is evaluated analytically using the displacement-discontinuity fundamental solutions (see Ref. [47]). The total perturbation displacement field is calculated by superposing the contributions from all fracture segments.

We discretize the governing equations in Eq. (1) by QFEM, QFEM-QPE, and MPSA, together with $\boldsymbol{u}_D$ provided by the reference displacement field and appropriate Neuman boundary conditions on $\Gamma_s^+$ and $\Gamma_s^-$ required in each benchmark problem. As a result, the corresponding displacement jump distribution along the fracture can be numerically recovered. This BEM approach to numerically implementing infinite-domain benchmark problems enables rigorous convergence analysis of SIFs in pure modes or mixed modes across different discretization methods. Fig. 3 depicts displacement fields and displacement jump distributions for Westergaard's mixed-mode problem computed by MPSA in the unit square domain, with a target mesh size of 0.025. The numerical results of Westergaard's mode-I, Westergaard's mode-II, and Sneddon's mode-I benchmarks by MPSA, QFEM, and QFEM-QPE are included in the code repository accompanying this paper [44].

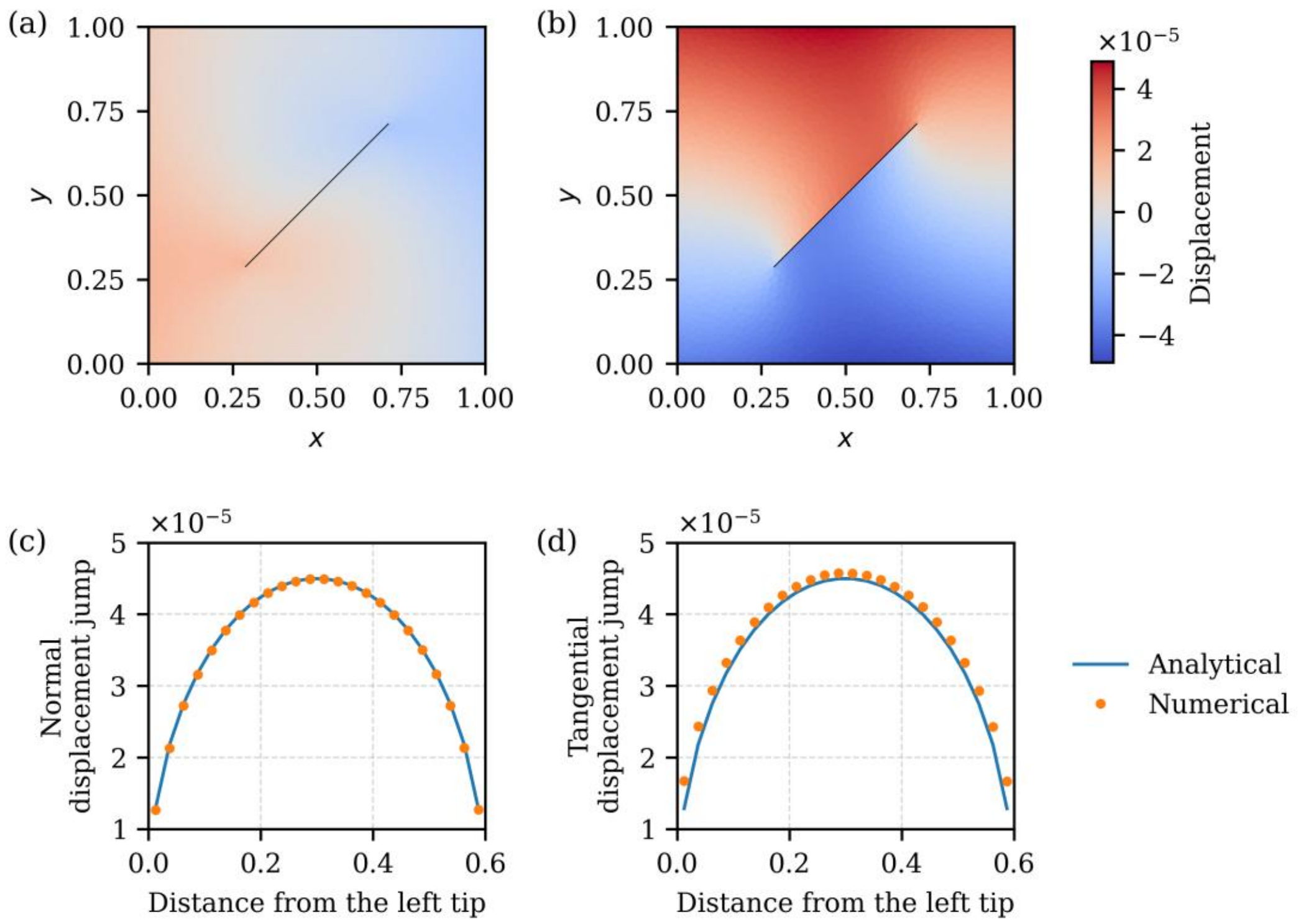

**Fig. 3.** Solution of displacement components, $u_x$ (a) and $u_y$ (b), in Westergaard's mixed-mode problems. Corresponding normal and tangential displacement jump profiles are shown in (c) and (d), respectively. MPSA discretization is adopted as an example.

### 4.6 Error estimation of discretization methods

We evaluate three types of errors from discretization methods, including displacement field errors over the matrix domain, displacement jump field errors over the fracture domain, and SIF estimation errors at fracture tips. First, let $\boldsymbol{u}_{num}$ and $\boldsymbol{u}_{ref}$ denote displacement fields on a given mesh interpolated by the numerical method and recovered by BEM, respectively. The relative $L^2$ error over the matrix domain $\Omega$ is defined as:

$$e_{\Omega}^{self} = \frac{\left(\int_{\Omega} \|\boldsymbol{u}_{num}-\boldsymbol{u}_{ref}\|^2 d\Omega\right)^{\frac{1}{2}}}{\left(\int_{\Omega} \|\boldsymbol{u}_{ref}\|^2 d\Omega\right)^{\frac{1}{2}}} \tag{37}$$

Second, let $\Gamma_s$ denote the fracture domain. The relative $L^2$ error between the numerical displacement jump field $[\![\boldsymbol{u}_{num}]\!]$ and the reference displacement jump field $[\![\boldsymbol{u}_{ref}]\!]$ over $\Gamma_s$ is defined as:

$$e_{\Gamma_s}^{self} = \frac{\left(\int_{\Gamma_s} \|[\![\boldsymbol{u}_{num}]\!]-[\![\boldsymbol{u}_{ref}]\!]\|^2 d\Gamma_s\right)^{\frac{1}{2}}}{\left(\int_{\Gamma_s} \|[\![\boldsymbol{u}_{ref}]\!]\|^2 d\Gamma_s\right)^{\frac{1}{2}}} \tag{38}$$

where $[\![\boldsymbol{u}_{ref}]\!]$ is obtained from the analytical solution of displacement jump in each benchmark. Third, the relative $L^2$ error is also used to evaluate the error of SIF estimates obtained from the scaling law. We calculate the discrete SIF, $K_{num,l}$, at tip $l$ from a discretization method and compare it with the reference SIF, $K_{ref}$ (see Eqs. (31) and (33)). The relative $L^2$ error of fracture tips in the whole domain, $e_K$, is defined as:

$$e_K = \frac{\left(\frac{1}{n_t}\sum_{l=1}^{n_t}\left(K_{num,l}-K_{ref}\right)^2\right)^{\frac{1}{2}}}{K_{ref}} \tag{39}$$

where $n_t$ is the number of fracture tips. For one straight fracture in the domain, $n_t$ is two.

## 5. Results and discussion

This section consists of three subsections. First, in Section 5.1, we perform global convergence analyses for QFEM, QFEM-QPE, and MPSA, followed by convergence analyses for extracted SIFs through the two-point displacement correlation method using the conventional nearest-tip sampling strategy. Then, Section 5.2 shows the performance of our proposed scaling law on SIF estimation in all three numerical discretization methods through four benchmark problems. Finally, in Section 5.3, we conduct convergence analyses of SIFs in these four benchmark problems with fitted parameters of the scaling law.

### 5.1 Failure of nearest-tip sampling

We evaluate convergence behavior of displacement solutions within the whole fractured domain across different discretization methods. Details of the convergence analysis are explained as follows. The same square computational domain containing a central horizontal fracture is used for all three discretization methods. The relative errors of displacement fields are calculated according to Eq. (37) by the numerical integration procedure over the matrix domain. The boundary element method introduced in Section 4.5.3 is used to compute the reference displacement at quadrature points in the given mesh for numerical integration. At the same quadrature points, the numerical displacement is interpolated directly by shape functions in QFEM and QFEM-QPE, while we reconstruct the displacement field in MPSA through volume-weighted averaging of cell-centered displacement at mesh nodes followed by the linear interpolation within each cell. We also compute the relative errors of displacement jumps in Eq. (38) through the numerical integration over the one-dimensional (1D) fracture domain. Due to the nonlinear mapping in QPEs, a $7 \times 7$ tensor-product Gaussian quadrature rule through the Duffy transformation is employed over the 2D triangular mesh for all three numerical discretization methods to ensure accurate approximation of integration [48, 49]. For the 1D fracture domain, the seven-point Gauss-Legendre quadrature rule is used for the numerical integration of displacement jumps in these three discretization methods [50]. Fig. 4 shows the convergence results for displacement fields over the matrix domain and for displacement jump fields along the fracture. The averaged convergence orders of displacement fields are estimated as 1.09, 1.48, and 1.21 for QFEM, QFEM-QPE, and MPSA, respectively, while the corresponding orders of displacement jump fields are 0.99, 1.34, and 0.97. For smooth problems, QFEM and MPSA are expected to exhibit third- and second-order convergence in the $L^2$ norm, respectively. However, in our fracture-embedded problems, the observed convergence rates of all three discretization methods are substantially degraded because the fracture-tip singularity lowers the regularity of the displacement field. In addition, we can see that the displacement jump fields from QFEM and MPSA both show first-order convergence. As shown in Fig.

4, the use of QPEs improves the convergence rates of QFEM for displacement fields and displacement jump fields, but it also introduces slight oscillations on the convergence curves. Compared with QFEM, MPSA shows higher convergence rates for the displacement field and features in weaker convergence degradation caused by fracture-tip singularity.

The convergence behavior of SIFs is presented in Fig. 5 through the two-point displacement correlation method using the conventional nearest-tip sampling strategy, with errors evaluated against the analytical SIF solutions for the corresponding benchmark problems (see Eqs. (31) and (33)). In this approach, SIFs are evaluated using the two sampled solutions of displacement jumps closest to each fracture tip. To account for the influence of fracture orientation on geometry and mesh, the fracture is rotated by the angle $\beta$, as illustrated in Fig. 6. SIFs are estimated using nearest-tip sampling for each of twelve fracture angles, with $\beta$ ranging from 0° to 165° in increments of 15°. As shown in Fig. 5, relative errors are reported as the mean over all considered fracture orientations with the shaded regions corresponding to $\pm 1$ standard deviation. Contrary to the convergent displacement fields shown in Fig. 4, the SIF errors do not decrease with mesh refinement for any of the three discretization methods. This failure of convergence is consistently observed in benchmark problems of mode I, mode II, and mixed modes. Compared with standard QFEM, QPEs reduce the magnitude of SIF errors in all benchmark cases. However, they do not produce convergent SIF predictions. This is because the $\sqrt{r}$-type behavior of displacement in the immediate vicinity of the fracture tip is not fully captured by numerical solutions from any of the tested discretization methods. Nearest-tip sampling consistently draws on this poorly resolved region and prevents accurate recovery of SIFs, as the sampled points move toward the tip with mesh refinement. The two-point displacement correlation method removes the first-order asymptotic truncation term for SIFs by extrapolation (see Eq. (9)). However, the extrapolation from the two-point method using inaccurate displacement samples near the fracture tip does not reduce the error in the SIF estimates and may instead propagate the displacement sampling errors into larger SIF estimation errors.

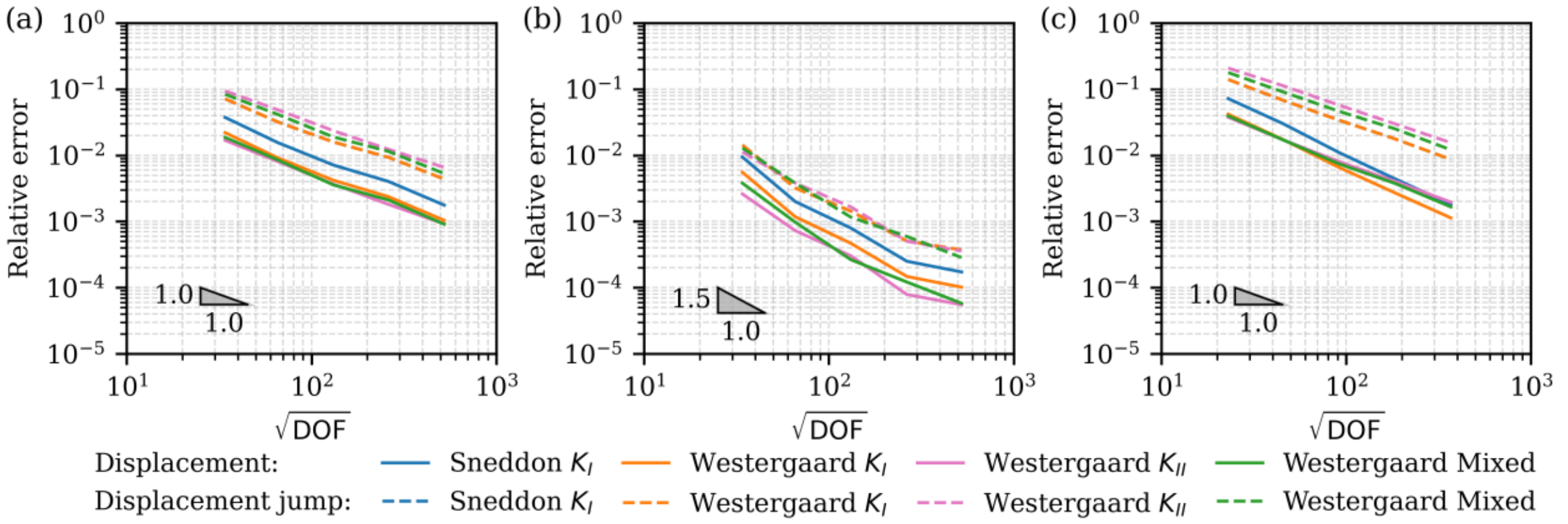

**Fig. 4.** Convergence of the numerical displacement and displacement jump fields for (a) QFEM, (b) QFEM-QPE, and (c) MPSA. Relative errors of displacement fields are evaluated over the matrix domain denoted by solid lines, while relative errors of displacement jump fields are calculated along the fracture denoted by dashed lines.

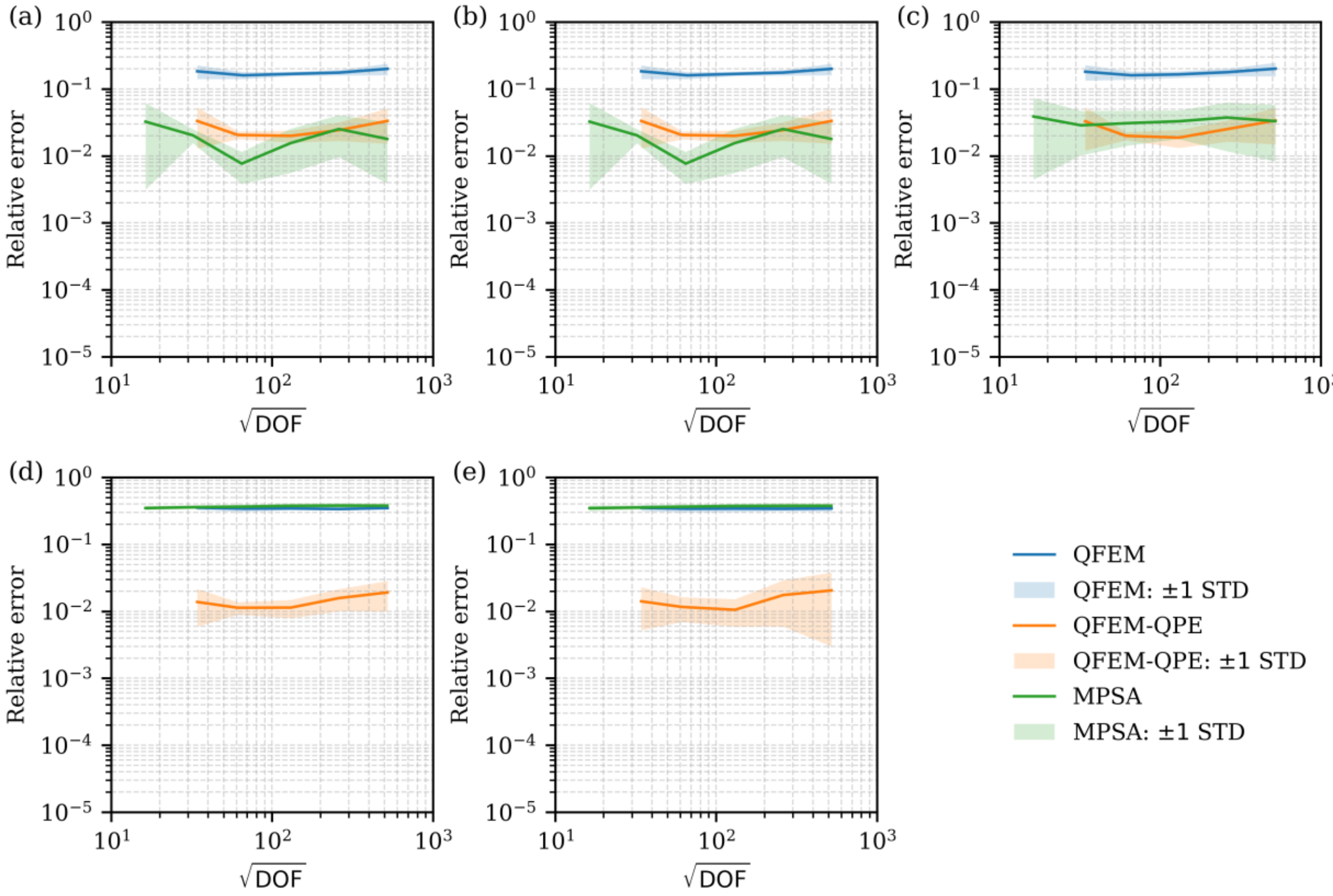

**Fig. 5.** Convergence behavior of SIF estimates using the nearest-tip sampling strategy across different benchmark cases and discretization methods. Panels (a)-(c) show the relative errors of the mode-I SIF for Sneddon's problem, Westergaard's mode-I problem, and Westergaard's mixed-mode problem,

respectively. Panels (d)-(e) show the relative errors of the mode-II SIF for Westergaard's mode-II problem and Westergaard's mixed-mode problem, respectively. Solid lines denote relative errors averaged over all fracture orientations, while shaded regions show ±1 standard deviation.

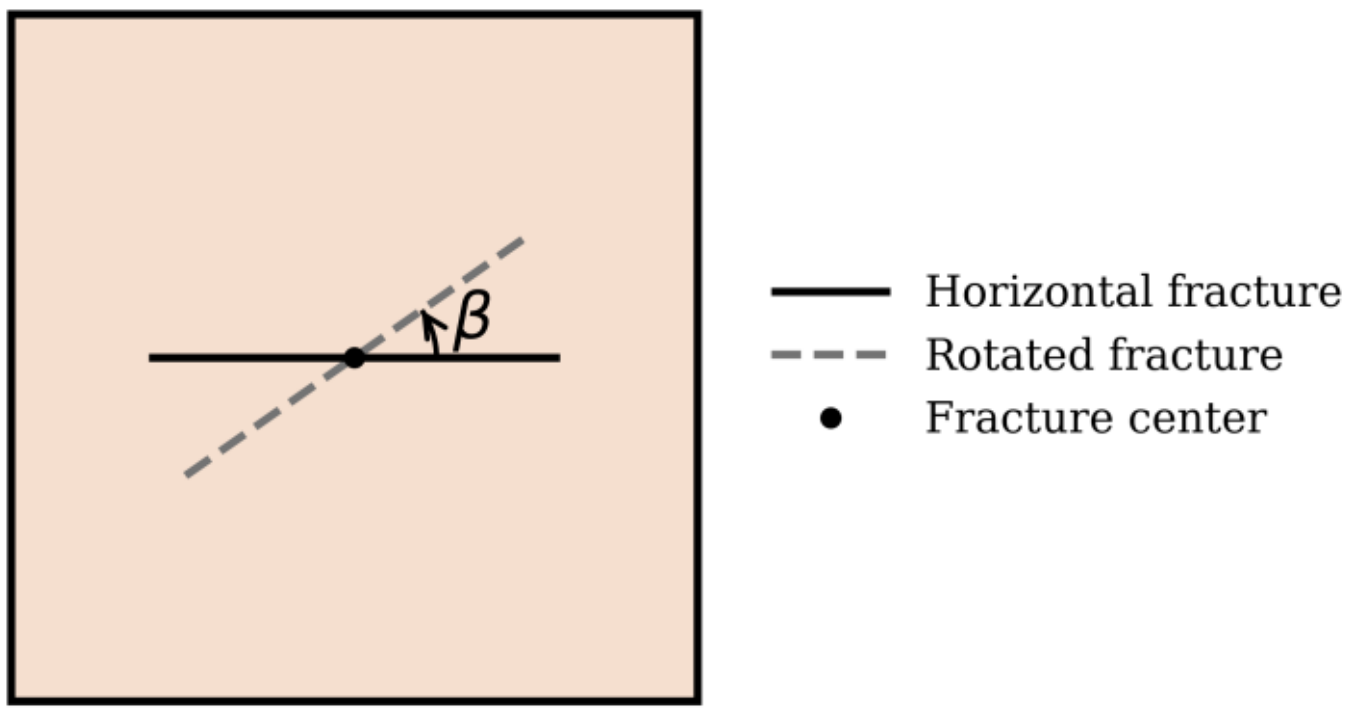

**Fig. 6.** Illustration of fracture orientation in the domain, with the fracture angle $\beta$ measured counterclockwise from the horizontal reference. The fracture center is fixed at the domain center.

### 5.2 Performance of the scaling law

Motivated by the failure of nearest-tip sampling, we evaluate how the proposed scaling law captures the mesh dependence of the optimal sampling distance for SIF estimation across different numerical discretization methods under mode-I, mode-II, mixed-mode loadings. The maps of relative SIF errors are constructed over a range of normalized sampling distance from the tip $(r/a)$ and normalized mesh size $(h/a)$ for each benchmark problem and discretization method. The geometry follows the centrally fractured domain with different fracture orientations, described in Section 5.1. For each pair $(r/a, h/a)$, the relative SIF errors obtained from all fracture orientations are averaged to produce a single error value. Figs. 7 and 8 show the calculated error maps for pure-mode and mixed-mode benchmark problems, respectively. Distinct low-error basins are found in all benchmark cases and discretization methods. A common feature observed across all error maps is that the low-error basins gradually move toward the fracture tip as the mesh is refined. This phenomenon demonstrates that the optimal sampling distance is mesh-size-dependent and cannot be represented by any fixed sampling locations.

Even though the low-error basins exist across all tested numerical methods, we find that their locations, widths, and error magnitudes are highly dependent on numerical discretization methods and fracture mode. Compared with MPSA, QFEM has much narrower low-error basins, with a strong sensitivity of the extracted SIFs to the optimal sampling distance. The introduction of QPEs in QFEM shifts the low-error basins closer to the fracture tip, suggesting an improved representation of the fracture-tip singularity. In addition, mode-II SIFs in MPSA have much broader low-error basins than mode-I SIFs, but with significantly higher error magnitudes.

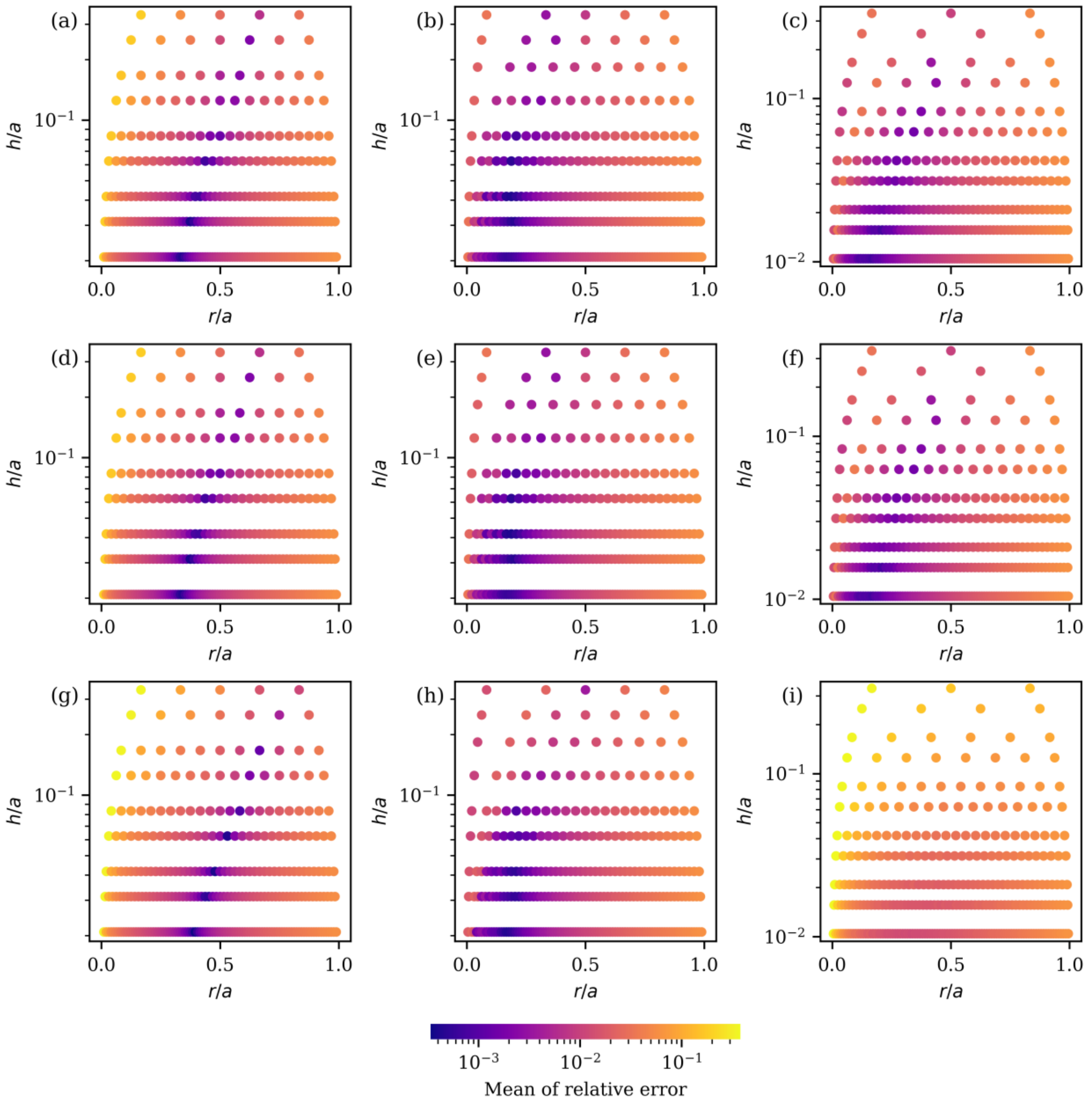

**Fig. 7.** The maps of relative SIF errors in the $(r/a, h/a)$ space for pure-mode benchmark problems using QFEM, QFEM-QPE, and MPSA. Errors are averaged over different fracture orientations at each point in the map. Panels (a)-(c) and panels (d)-(f) show relative errors for the mode-I SIF in Sneddon's problem and Westergaard's mode-I problem, respectively. Panels (g)-(i) show relative errors for the mode-II SIF in Westergaard's mode-II problem. Columns correspond to QFEM, QFEM-QPE, and MPSA, respectively.

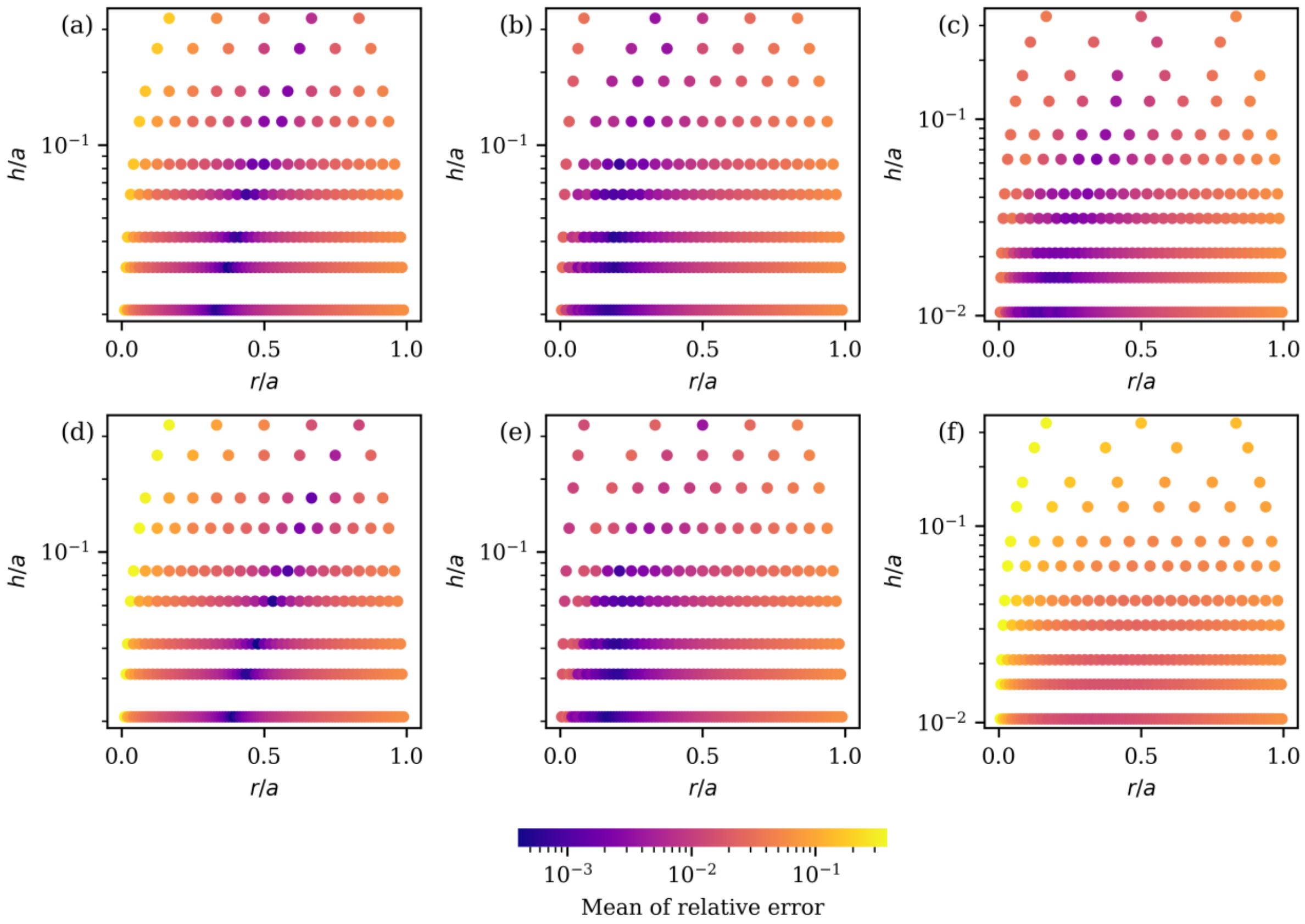

**Fig. 8.** The maps of relative SIF errors in the $(r/a, h/a)$ space for Westergaard's mixed-mode problem using QFEM, QFEM-QPE, and MPSA. Errors are averaged over different fracture orientations at each point in the map. Panels (a)-(c) and panels (d)-(f) show relative SIF error for mode I and mode II, respectively. Columns correspond to QFEM, QFEM-QPE, and MPSA, respectively.

To quantify the shift trend of optimal sampling distances, we perform two fitting procedures. First, we locate the minimized error at each mesh refinement level in the averaged error map. These optimal sampling locations are fitted using the proposed scaling law in Eq. (17), as shown in Figs. 9 and 10. The resulting data and fitting lines are shown as blue dots and blue dashed lines, respectively. The trajectory of the minimum-error location determined by the averaged error map is regarded as the global result. To further examine the robustness of the scaling law with respect to fracture orientation, we perform a second fitting procedure. For each fracture orientation, the optimal sampling location is independently determined at each given mesh. This produces a set of orientation-dependent optimal sampling locations at each mesh refinement level. We calculate the mean values and corresponding standard deviations of these optimal locations in Figs. 9 and 10, denoted as orange dots and orange shaded regions, respectively. The mean values are fitted using the scaling law (see orange dashed lines). This second fitting procedure evaluates how the orientation-dependent variability affects the fitted scaling-law parameters.

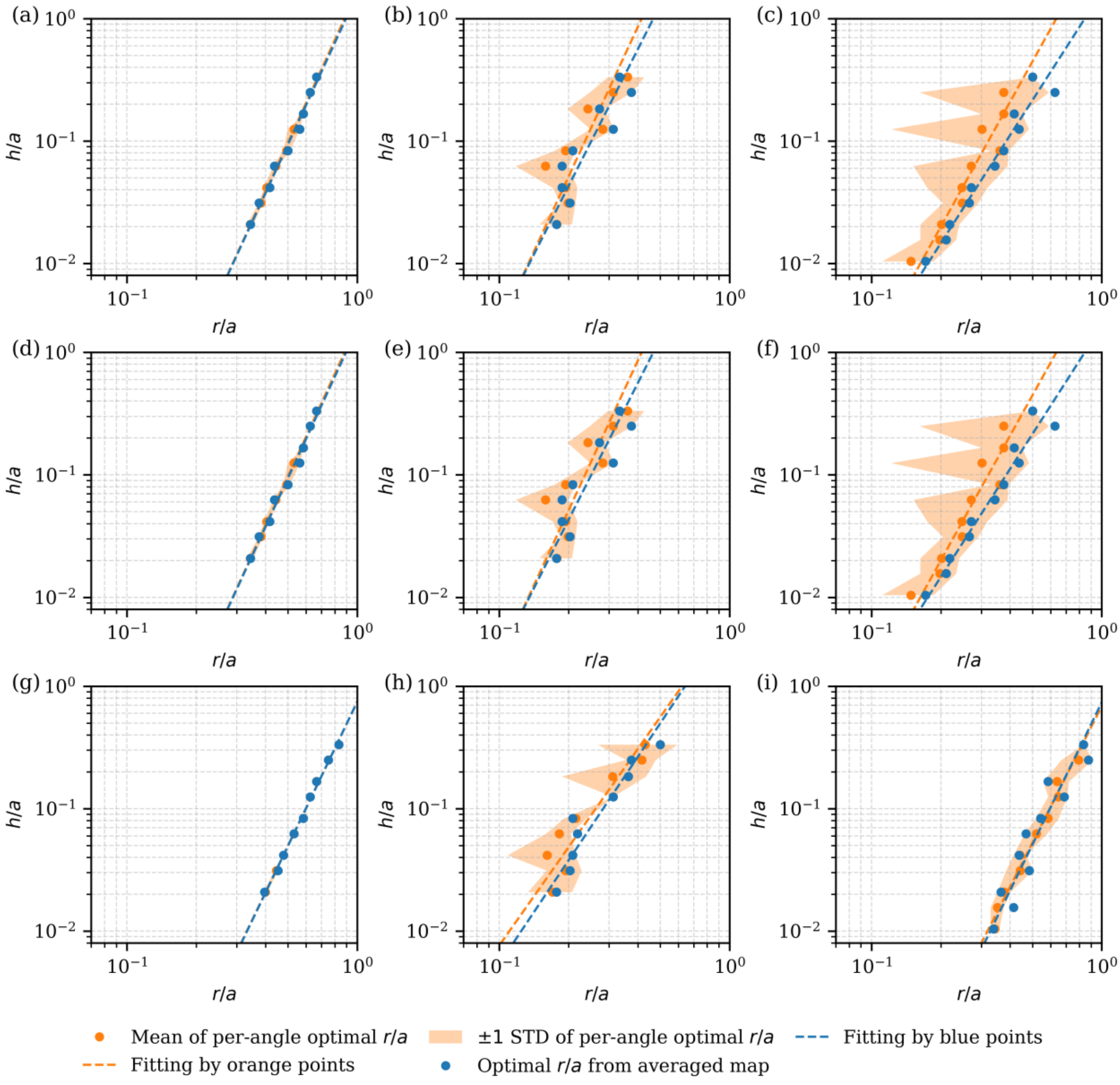

**Fig. 9.** The fitting performance of the scaling law relating the optimal sampling distance to the normalized mesh size for QFEM, QFEM-QPE, and MPSA in pure-mode benchmark problems. Panels (a)-(c) and panels (d)-(f) show fitting results of the mode-I SIF in Sneddon's problem and Westergaard's mode-I problem, respectively. Panels (g)-(i) show fitting results of the mode-II SIF in Westergaard's mode-II problem. Columns correspond to QFEM, QFEM-QPE, and MPSA, respectively.

The two fitting procedures have complementary purposes. The former fitting based on the averaged error maps characterizes the global minimum-error sampling trajectories in the $(r/a, h/a)$ space, while the latter fitting based on orientation statistics quantifies the sensitivity of the optimal sampling distance to fracture rotations. The numerical results in Figs. 9 and 10 show a very close agreement between these two fitting procedures. This confirms that our proposed scaling law captures the dominant trends of optimal sampling locations, whose fitted parameters should be largely independent of fracture orientations. However, in practice, arbitrary fracture orientations introduce variations in the geometry and the resulting mesh. Therefore, the standard deviations quantify the variability of the optimal sampling distance caused by fracture orientation variations for each discretization method. QFEM has the smallest standard deviations among these three discretization methods, as shown in Figs. 9 and 10. This indicates that QFEM has the most consistent minimum-error locations along the $r/a$ axis across different fracture orientations. The variability of the optimal sampling distance is also affected by the width of the low-error basins shown in Figs 7. and 8. With a broad and flat basin, small changes in the error profile may substantially shift the minimum-error location along the $r/a$ axis at a given mesh resolution. However, this shift may not change the SIF error values significantly. Therefore, broader basins provide higher tolerance to deviations from the optimal sampling distance, which translates into higher variability in the minimum-error locations.

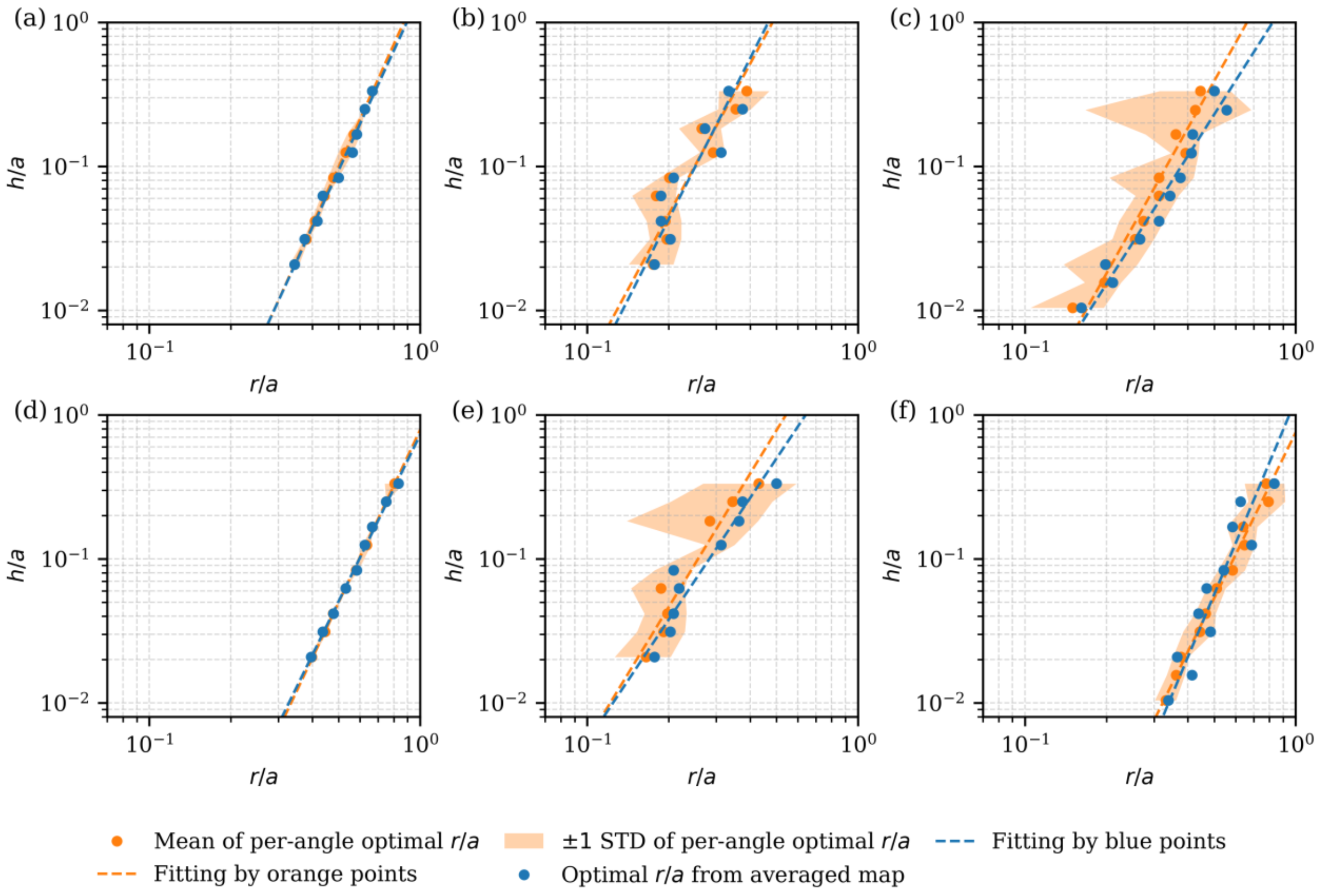

**Fig. 10.** The fitting performance of the scaling law relating the optimal sampling distance to the normalized mesh size for QFEM, QFEM-QPE, and MPSA in Westergaard's mixed-mode problem. Panels (a)-(c) and panels (d)-(f) show fitting results of SIFs for mode I and mode II, respectively. Columns correspond to QFEM, QFEM-QPE, and MPSA, respectively.

**Table 1.** Fitted parameters of the scaling law over averaged error maps for mode-I SIFs in QFEM, QFEM-QPE, and MPSA.

| **Benchmarks** | **QFEM** | | **QFEM-QPE** | | **MPSA** | |
|---|---|---|---|---|---|---|
| | ***C*** | ***α*** | ***C*** | ***α*** | ***C*** | ***α*** |
| **Sneddon's problem** | 0.894 | 0.245 | 0.467 | 0.269 | 0.844 | 0.341 |
| **Westergaard's mode-I problem** | 0.894 | 0.245 | 0.467 | 0.269 | 0.844 | 0.341 |
| **Westergaard's mixed-mode problem** | 0.894 | 0.245 | 0.467 | 0.269 | 0.823 | 0.338 |

**Table 2.** Fitted parameters of the scaling law over averaged error maps for mode-II SIFs in QFEM, QFEM-QPE, and MPSA.

| Benchmarks | QFEM | | QFEM-QPE | | MPSA | |
|---|---|---|---|---|---|---|
| | $C$ | $\alpha$ | $C$ | $\alpha$ | $C$ | $\alpha$ |
| **Westergaard's mode-II problem** | 1.079 | 0.255 | 0.642 | 0.356 | 1.082 | 0.259 |
| **Westergaard's mixed-mode problem** | 1.087 | 0.260 | 0.642 | 0.356 | 0.948 | 0.223 |

Tables 1 and 2 show the fitted parameters of the scaling law over averaged error maps for mode-I and mode-II SIFs, respectively. A notable finding is that for a given numerical discretization method the fitted parameters remain similar across all different benchmark problems within the same fracture mode. This numerical observation suggests that the scaling law is primarily governed by the fracture mode, instead of the specific loading configurations from different benchmark cases for the same mode. Furthermore, according to these fitted parameters, the scaling law is both mode-dependent and method-dependent. Among the three discretization methods, QFEM-QPE has the smallest scaling constant $C$ in both fracture modes. This is consistent with the role of QPEs in improving representation of fracture-tip singularity. The scaling exponent $\alpha$ is related to the convergence rates of discretization methods degraded due to asymptotic truncation error and numerical discretization error. Its variation does not follow a common trend between the two fracture modes: MPSA has the largest $\alpha$ in mode I, while the largest $\alpha$ in mode II is from QFEM-QPE.

### 5.3 Convergence of SIF estimates

We investigate the convergence behavior of SIF estimates in QFEM, QFEM-QPE, and MPSA using the scaling law in Eq. (17) with the fitted parameters in Tables 1 and 2. In contrast to the non-convergent results by the nearest-tip sampling strategy in Fig. 5, the scaling law enables a convergent SIF computation with mesh refinement (see Fig. 11). This improvement is validated in both mode I and mode II across the three discretization methods. For mode-I estimates, QFEM, QFEM-QPE, and MPSA have similar convergence performance. In mode II, QFEM and QFEM-QPE have lower magnitudes of relative errors than MPSA does, while MPSA exhibits a very stable convergence rate of SIF estimates.

The scaling law restores the convergent SIF estimates across all three discretization methods. However, the calculated convergence curves are not perfectly linear in the log-log scale. Different from the displacement fields, SIFs are derived quantities calculated from the post-processing procedure based on the numerical solution. Additional sources of variations that are negligible in the displacement convergence analysis are inherited and amplified when SIF estimates are extracted. Consequently, the convergence behavior of relative SIF errors is not only affected by convergence rates of displacement solutions but is also constrained by the complicated interactions among numerical discretization errors, asymptotic truncation errors, discrete sampling locations projected by the scaling law, and meshing quality variations induced by fracture orientations.

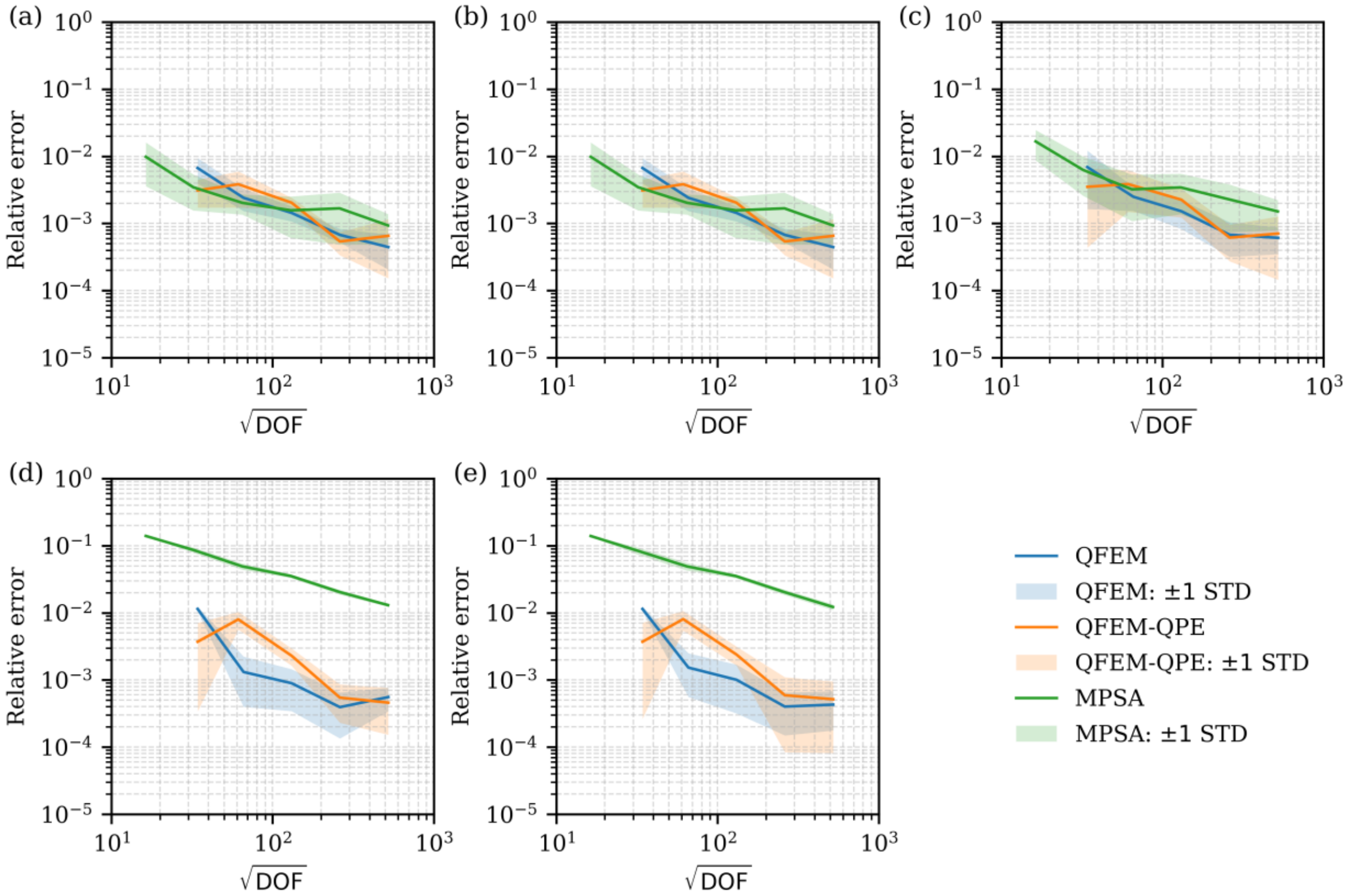

**Fig. 11.** The convergence behavior of mode-I SIF estimates in panels (a)-(c) and of mode-II SIF estimates in panels (d)-(e) across different discretization methods using scaling-law sampling. Panels (a), (b), and (d) are calculated from Sneddon's problem, Westergaard's mode-I problem, and Westergaard's mode-II problem, respectively. Panels (c) and (e) are both calculated from Westergaard's mixed-mode problem.

## 6. Conclusions

We carried out a numerical study for fracture SIF computations using two-point displacement correlation methods in three numerical discretization methods, QFEM, QFEM-QPE, and MPSA. Numerical simulations were performed for pure-mode and mixed-mode benchmark cases for fracture opening and sliding including Westergaard-type traction-free fracture problems and Sneddon pressurized fracture problems. First, our numerical experiments demonstrate that the conventional nearest-tip sampling strategy does not guarantee the convergence of SIF estimates for any of the three discretization methods. Then, based on a balance between the numerical discretization error and the asymptotic truncation error, we proposed a scaling law as a new sampling strategy along the fracture towards the fracture tip for SIF computations. The maps of relative SIF errors were constructed in terms of the normalized sampling distance from the tip and the normalized mesh size. The fitted parameters of the scaling law depend on fracture modes and discretization methods, but not on the particular benchmark problems. Finally, the numerical results indicate that the proposed scaling law enables convergent SIF estimates in both pure and mixed modes across all three discretization methods.

## Appendix A. Asymptotic truncation error of the two-point displacement correlation method

For either fracture mode $\varphi \in \{I, II\}$, the exact SIF, $K_\varphi$ , can be evaluated by Eq. (6) at the sampling distance $r_i$ with the single-point estimate $K_\varphi^{sp}(r_i)$ from Eq. (8) and the fourth-order truncation, yielding

$$K_\varphi = K_\varphi^{sp}(r_i) + g_1^\varphi r_i + g_2^\varphi r_i^2 + g_3^\varphi r_i^3 + O(r_i^4). \quad \text{(A.1)}$$

Substituting Eq. (A.1) at $r_i$ and $r_j$ into Eq. (9), the two-point displacement correlation method gives

$$K_\varphi^{tp}(r_i, r_j) = K_\varphi + g_2^\varphi r_i r_j + O\left(r_i r_j (r_i + r_j)\right). \quad \text{(A.2)}$$

Eq. (A.2) indicates that the leading truncation errors in two-point estimates are proportional to the product $r_i r_j$. We define $r = \max\{r_i, r_j\}$ to characterize the two sampling locations by a single length scale. Without loss of generality, let $r_i = r$ and $r_j = \rho r$, where $\rho$ is the ratio determined by the available sampling location and $0 < \rho < 1$. In this way, Eq. (A.2) is rewritten as

$$K_\varphi^{tp} = K_\varphi + g_2^\varphi \rho r^2 + O(\rho(1 + \rho) r^3). \quad \text{(A.3)}$$

Since $\rho$ remains bounded between 0 and 1, the asymptotic truncation error $e_{\varphi}^{asym}$ for either fracture mode satisfies

$$e_{\varphi}^{asym} = \left|K_{\varphi}^{tp} - K_{\varphi}\right| = O(r^2). \tag{A.4}$$

**Appendix B. Numerical discretization error of the two-point displacement correlation method**

Analogous to Eq. (8) evaluated by exact displacement $\boldsymbol{u}(r_i)$, the single-point SIF estimates for mode I and mode II, $K_{I,h}^{sp}$ and $K_{II,h}^{sp}$, evaluated by the numerical displacement field $\boldsymbol{u}_h(r_i)$ at $r_i$ are given by

$$\begin{aligned} K_{I,h}^{sp}(r_i) &= \frac{\mu}{(\kappa+1)}\sqrt{\frac{2\pi}{r_i}}\,[\![\boldsymbol{u}_h(r_i)]\!]_n, \\ K_{II,h}^{sp}(r_i) &= \frac{\mu}{(\kappa+1)}\sqrt{\frac{2\pi}{r_i}}\,[\![\boldsymbol{u}_h(r_i)]\!]_t. \end{aligned} \tag{B.1}$$

For either fracture mode $\varphi \in \{I, II\}$, the signed discretization error of the single-point SIF estimate at $r_i$, $e_{\varphi}^{sp}$, is defined as

$$e_{\varphi}^{sp}(r_i) = K_{\varphi,h}^{sp}(r_i) - K_{\varphi}^{sp}(r_i). \tag{B.2}$$

Using $\eta_{\varphi}^{num}(r)$ defined in Eq. (12), Eq. (B.2) can be rewritten as

$$e_{\varphi}^{sp}(r_i) = \frac{\mu}{(\kappa+1)}\sqrt{\frac{2\pi}{r_i}}\,\eta_{\varphi}^{num}(r_i). \tag{B.3}$$

The SIF estimate, $K_{\varphi,h}^{tp}$, from the two-point displacement correlation method by the numerical displacement field is given by

$$K_{\varphi,h}^{tp}(r_i, r_j) = K_{\varphi,h}^{sp}(r_i) + \frac{r_i}{r_j - r_i}\left(K_{\varphi,h}^{sp}(r_i) - K_{\varphi,h}^{sp}(r_j)\right). \tag{B.4}$$

The discretization error of two-point SIF estimates, $e_{\varphi}^{num}$, is defined as

$$e_{\varphi}^{num} = \left|K_{\varphi,h}^{tp} - K_{\varphi}^{tp}\right|. \tag{B.5}$$

Substituting Eqs. (B.2), (B.4) and (9) into Eq. (B.5), $e_{\varphi}^{num}$ is therefore rewritten as

$$e_{\varphi}^{num}(r_i, r_j) = \left|e_{\varphi}^{sp}(r_i) + \frac{r_i}{r_j - r_i}\left(e_{\varphi}^{sp}(r_i) - e_{\varphi}^{sp}(r_j)\right)\right|. \tag{B.6}$$

Using Eq. (B.3), Eq. (B.6) can be expressed in terms of $\eta_{\varphi}^{num}$ as

$$e_{\varphi}^{num}(r_i, r_j) = \left|\frac{\mu}{(\kappa+1)}\frac{\sqrt{2\pi}}{r_j - r_i}\left(\frac{r_j\eta_{\varphi}^{num}(r_i)}{\sqrt{r_i}} - \frac{r_i\eta_{\varphi}^{num}(r_j)}{\sqrt{r_j}}\right)\right|. \tag{B.7}$$

Following Appendix A for two sampling locations by a single length scale, we assume $r_i = r$ and $r_j = \rho r$, where $0 < \rho < 1$. Eq. (B.7) is rewritten as

$$e_{\varphi}^{num}(r_i, r_j) = \frac{\mu\sqrt{2\pi}}{(\kappa+1)} Q(\rho)\left|c_{\varphi}\right| h^s r^{-0.5}, \quad \text{(B.8)}$$

where

$$Q(\rho) = \left(\frac{\rho - \rho^{-0.5}}{\rho - 1}\right). \quad \text{(B.9)}$$

With the assumption in Eq. (13), $Q(\rho)$ remains finite when $\rho$ approaches one. Finally, the numerical discretization error $e_{\varphi}^{num}$ for either fracture mode is predicted to scale as

$$e_{\varphi}^{num} = O(h^s r^{-0.5}). \quad \text{(B.10)}$$

## CRediT authorship contribution statement

**Ronghao Cui:** Conceptualization, Methodology, Software, Validation, Investigation, Writing - Original Draft. **Jan Martin Nordbotten:** Conceptualization, Supervision, Methodology, Writing - review & editing. **Ivar Stefansson:** Conceptualization, Supervision, Methodology, Software, Writing - review & editing. **Inga Berre:** Conceptualization, Supervision, Methodology, Project administration, Funding acquisition, Writing - review & editing.

## Declaration of competing interest

The authors declare that they have no known competing financial interests or personal relationships that could have appeared to influence the work reported in this paper.

## Declaration of generative AI and AI-assisted technologies in the manuscript preparation process

During the preparation of this work, the authors used the OpenAI-GPT5.6 model to generate editing suggestions. The use of this tool was limited to improving grammar, clarity, and readability of the manuscript. After using this tool, the authors reviewed and edited the content as needed and take full responsibility for the content of the published article.

## Data availability

The data and source code used to generate the results presented in this work are available as a Docker image on Zenodo, as documented by Cui et al. [44].

## Acknowledgements

This result is part of a project that has received funding from the European Research Council (ERC) under the European Union's Horizon 2020 research and innovation programme (grant agreement No 101002507).